\documentclass[leqno,a4paper]{article}

\usepackage{amssymb,amsmath,amsthm}

\title{\textbf{Differential }$p$\textbf{-forms and }$q$\textbf{-vector fields }
\textbf{with constant coefficients}}
\author{\textsc{J. Mu\~{n}oz Masqu\'e}
\and \textsc{L. M. Pozo Coronado}
\and \textsc{M. E. Rosado Mar\'{\i}a}}

\date{}

\newtheorem{theorem}{Theorem}[section]
\newtheorem{proposition}[theorem]{Proposition}
\newtheorem{lemma}[theorem]{Lemma}
\newtheorem{corollary}[theorem]{Corollary}

\theoremstyle{remark}
\newtheorem{remark}[theorem]{Remark}
\newtheorem{notation}[theorem]{Notation}
\newtheorem{definition}[theorem]{Definition}
\newtheorem{example}[theorem]{Example}

\begin{document}

\maketitle

\begin{abstract}
\noindent Differential $p$-forms and $q$-vector fields with constant
coefficients are studied. Differential $p$-forms of degrees $p=1,2,n-1,n$
with constant coefficients on a smooth $n$-dimensional manifold $M$ are
characterized. In the contravariant case, the obstruction for a $q$-vector
field $V_q$ to have constant coefficients is proved to be the Schouten-Nijenhuis
bracket of $V_q$ with itself. The $q$-vector fields with constant
coefficients of degrees $q=1,2,n-1,n$ are also characterized.
The notions of differential $p$-forms and $q$-vector fields
with conformal constant coefficients are introduced.
For arbitrary degrees $p$ and $q$, such
differential $p$-forms and $q$-vector fields are seen to be the solutions
to two second-order partial differential systems on $J^2(M,\mathbb{R}^n)$,
which are reducible to two first-order partial differential systems
by adding variables. Computational aspects in solving
these systems are discussed and examples and applications are also given.
\end{abstract}


\medskip


\section{Definitions and examples}
The notion of ``having constant coefficients'' plays an important role
in both Geometry and Analysis; for example, when does a differential operator
have constant coefficients? (\cite{GuS}), or when an exterior differential
system (\cite[Ch.\ XVIII, \S 9]{Dieu}) is defined by forms with constant coefficients,
etc. In Geometry this notion is usually related to concepts such as locally
symmetric spaces, flat semi-Riemannian manifolds, etc. Below we tackle this notion
in the case of $p$-differential forms and $q$-vector fields.
First of all we introduce the formal definitions and show some examples.
\subsection{The covariant case}
\begin{definition}
\label{def1}
Let $M$ be a smooth manifold of dimension $n$. A differential form
$\omega _p\in \Omega ^p(M)$ is said to have \emph{constant coefficients} on
a neighbourhood of the point $x\in M$ if there exists a system of coordinates
$(u^1,\dotsc,u^n)$ centred at $x$, such that all the functions
$\omega _p\left(
\tfrac{\partial }{\partial u^{i_1}},\dotsc,\tfrac{\partial }{\partial u^{i_p}}
\right) $,
$1\leq i_1<\ldots <i_p\leq n$, are constant.
\end{definition}
\begin{notation}
\label{not1}
For every $p\in \mathbb{N}^+$, let $\mathcal{I}_p^n$ denote
the set of all multi-indices $I=(i_1,\dotsc,i_p)\in \mathbb{N}^p$ such
that $1\leq i_1<\ldots<i_p\leq n$. If $(u^1,\dotsc,u^n)$ is a system
of $C^{\infty} $ functions on $M$, we write
$du^I=du^{i_1}\wedge \cdots\wedge du^{i_p}$,
$\forall I\in \mathcal{I}_p^n$.
\end{notation}
\begin{proposition}
\label{propos-1}
Every form $\omega _p$\ with constant coefficients is closed.
In this case, the differential system
$D=\{X\in T_xM:i_{X}(\omega _p)_x=0,x\in M\} $ has a locally
constant rank and it is involutive.
\end{proposition}
\begin{proof}
If $\omega _p$ has constant coefficients, by writing $\omega _p$ on a
coordinate system as in Definition \ref{def1}, we have $\omega _p
=\sum\nolimits_{I\in \mathcal{I}_p^n}\lambda_Idu^I$,
$\lambda _I\in \mathbb{R}$, and by taking its exterior differential,
we have $d\omega _p=\sum_{I\in \mathcal{I}_p^n}d\lambda _I\wedge du^I=0$.

The tangent vectors $X\in T_xM$\ in $D$ are in one-to-one correspondence
with the solutions
\[
\mathcal{X}=\{ (X^1,\dotsc,X^n)\in \mathbb{R}^n\} ,\quad X^i
=du^i(X),\quad 1\leq i\leq n,
\]
to the following system of $\binom{n}{p-1}$ linear equations
with constant coefficients:
\[
0=\sum\nolimits_{I\in \mathcal{I}_p^n}(-1)^{h-1}X^h\lambda _I
\left( du^{i_1}\right) _x\wedge \cdots \wedge
\widehat{\left( du^{i_h}\right) _x}\wedge \cdots \wedge
\left( du^{i_p}\right) _x.
\]
As the vector space $\mathcal{X}$ is independent of the point $x$,
it follows that $D$ has a basis of vector fields with constant coefficients
$V^1,\dotsc,V^k$ on the domain of the system $(u^1,\dotsc,u^n)$. Hence
$[V^i,V^j]=0$, $\forall i,j=1,\dotsc,k$, and accordingly, $D$ is involutive.
\end{proof}
\begin{remark}
\label{remark0}
In degrees $1$, $2$, $n-1$ and $n$ the property of having
constant coefficients is generic in the space of closed forms, as proved in
the next four examples.
\end{remark}
\begin{example}
\label{example1}
A $1$-form $\omega _1$ that does not vanish at $x$ has
constant coefficients if and only if $\omega _1$ is closed. If $\omega_1$
has constant coefficients, then it is closed as follows from Proposition
\ref{propos-1}; hence $\omega _1$ is locally exact: $\omega _1=df$. We can
suppose $f(x)=0$ and $\frac{\partial f}{\partial x^1}(x)\neq 0$, permuting
indices if necessary. In this case the functions $u^1=f$, $u^h=x^h$,
$2\leq h\leq n$, are a coordinate system centred at the point $x$, as the
Jacobian of the system $(u^i)_{i=1}^n$ with respect to $(x^j)_{j=1}^n$
at $x$ is $\tfrac{\partial(u^1,\dotsc,u^n)}{\partial(x^1,\dotsc,x^n)}(x)
=\tfrac{\partial f}{\partial x^1}(x)$.
\end{example}
\begin{example}
\label{example2}
A form of maximum degree $\omega _n\in \Omega ^n(M)$ has
constant coefficients on a neighbourhood of $x\in M$ if, either
$\omega _n$ vanishes on a neighbourhood of $x$, or $\omega _n$
does not vanish at $x$. If $\omega _n=fdx^1\wedge \cdots \wedge dx^n$
and $f(x)\neq 0$, then by defining $u^1=\int _0^{x^1}
f\left( t,x^2,\dotsc,x^n\right) dt$, $u^i=x^i$, $2\leq i\leq n$,
we have $\omega _n=du^1\wedge \cdots \wedge du^n$.
\end{example}
\begin{example}
\label{example3}
A $2$-form $\omega _2$ of constant class
(cf.\ \cite[Appendix 4, \S3.5]{Libermann-Marle}) on a neighbourhood
of $x\in M$ has constant coefficients if and only it is closed,
as follows from Darboux's theorem (cf.\ \cite[VI, \S 4.5]{Godbillon}).
\end{example}
\begin{example}
\label{example4}
A $(n-1)$-form $\omega _{n-1}\in \Omega ^{n-1}(M)$ that does not
vanish at $x\in M$ has constant coefficients if and only it is closed. In
fact, if $\omega _n$ is a volume form on an open neighbourhood $U$ of
$x\in M$, then there exists a vector field $X\in \mathfrak{X}(U)$ such that
$i_X\omega _n=\omega _{n-1}$, and $X_x\neq 0$; hence, there exist
coordinates $(x^i)$ centred at $x$ reducing $X$ to a canonical form, i.e.,
$X=\frac{\partial }{\partial x^1}$. On these coordinates:
$\omega _n=\rho dx^1\wedge\ldots\wedge dx^n$ and
$\omega_{n-1}=\rho dx^2\wedge \cdots\wedge dx^n$. By differentiating
we thus have
$0=d\omega_{n-1}
=\frac{\partial\rho }{\partial x^1}dx^1\wedge \cdots \wedge dx^n$.
Accordingly, $\rho $ depends only on $x^2,\dotsc,x^n$,
and we can conclude by virtue of Example \ref{example2}.
\end{example}
\begin{remark}
\label{remark1}
The cases not included in the four examples above are $n\geq 5$,
$3\leq p\leq n-2$.
\end{remark}
\begin{example}
\label{example4.5}
A differential $p$-form $\omega _p$ on $M$ of constant class
equal to $p$
has constant coefficients, as every $x\in M$ admits
a coordinate neighbourhood $(U;x^1,\dotsc,x^n)$ such that
$\omega _p|_U=dx^1\wedge \cdots \wedge dx^p$; see \cite[VI, 3.4]{Godbillon}.
\end{example}
\subsection{The contravariant case}
\begin{definition}
\label{def0}
Let $M$ be a smooth manifold of dimension $n$. A $q$-vector field
$V_q\in \wedge ^q\mathfrak{X}(M)$ is said to have \emph{constant}
\emph{coefficients} on a neighbourhood of the point $x\in M$ if there exists a
system of coordinates $(u^1,\dotsc,u^n)$ centred at $x$, such that all the
functions $V_q\left( du^{i_1},\dotsc,du^{i_q}\right) $,
$1\leq i_1<\ldots <i_q\leq n$, are constant.
\end{definition}
\begin{notation}
\label{not2}
We set $\partial u^I
=\frac{\partial }{\partial u^{i_1}}\wedge \cdots \wedge
\frac{\partial }{\partial u^{i_q}}$, $\forall I\in \mathcal{I}_q^n$.
\end{notation}
\begin{example}
\label{example5}
A $1$-vector field $X\in \mathfrak{X}(M)$ that does not vanish
at $x$ has constant coefficients on a neighbourhood of this point, as there
exists a system of coordinates $(U;x^1,\dotsc,x^n)$ centred at $x$, such
that: $X|_U=\frac{\partial }{\partial x^1}$.
\end{example}
\begin{proposition}
\label{propos0}
If a $q$-vector field $V_q$ on $M$ has constant coefficients,
then Schouten-Nijenhuis bracket $V_q$ with itself vanishes, i.e.,
$[V_q,V_q]=0$. In this case, the Pfaffian system
$P=\{w \in T_x^\ast M:i_w(V_q)_x=0,x\in M\} $ has a locally
constant rank and for every $x\in M$ there exists an open
neighbourhood $U_x$ such that $\Gamma (U_x,P)$ admits a basis
of closed $1$-forms; hence, $P$ is integrable.
\end{proposition}
\begin{proof}
If $V_q=\sum _{I\in \mathcal{I}_q^n}\lambda _I\partial u^I$,
$\lambda _I\in \mathbb{R}$, in a coordinate system $(U;u^1,\dotsc,u^n)$,
the first part of the statement follows by virtue of the following formula
(cf.\ \cite{Michor}):
\[
\begin{array}
[c]{rl}
\lbrack A,B]= &
{\textstyle\sum_{i,j}}
(-1)^{i+j}[X_i,Y_j]\wedge X_1\wedge \cdots \wedge \widehat{X_i}
\wedge \cdots \wedge X_q\wedge Y_1\wedge\\
& \multicolumn{1}{r}{\cdots \wedge \widehat{Y_j}\wedge \cdots \wedge Y_r,}\\
A= & X_1\wedge \cdots \wedge X_q,\\
B= & Y_1\wedge \cdots \wedge Y_r,\\
\forall X_i,Y_j\in \mathfrak{X}(M), & 1\leq i\leq q,1\leq j\leq r.
\end{array}
\]
The $1$-covectors $w\in T_x^\ast M$ in $P$ are in one-to-one correspondence
with the solutions
\[
\mathcal{W}=\{ (w^1,\dotsc,w^n)\in \mathbb{R}^n\} ,
\quad w^i=w\left(
\tfrac{\partial }{\partial u^i}\right) _x,\quad 1\leq i\leq n,
\]
to the following system of $\binom{n}{q-1}$ linear equations
with constant coefficients:
\[
0=\sum\nolimits_{I\in \mathcal{I}_q^n}(-1)^{h-1}w^h\lambda _I
\left(
\tfrac{\partial }{\partial u^{i_1}}\right) _x\wedge \cdots \wedge
\widehat{\left( \tfrac{\partial }{\partial u^{i_h}}\right) _x}
\wedge \cdots \wedge\left( \tfrac{\partial }{\partial u^{i_q}}\right) _x.
\]
As the vector space $\mathcal{W}$ is independent of the point $x$, it follows
that $P$\ has a basis of $1$-forms $\zeta ^1,\dotsc,\zeta ^k$ with constant
coefficients and consequently $d\zeta^j=0$, $1\leq j\leq k$. Hence
\[
d\left( \sum\nolimits_{j=1}^kf_j\zeta ^j\right)
=\sum\nolimits_{j=1}^kdf_j\wedge \zeta ^j,
\quad \forall f_j\in C^\infty (M),\;1\leq j\leq k.
\]
This proves that $P$ is integrable.
\end{proof}
\begin{corollary}
\label{cor}
If $n=\dim M$ is even, then a $2$-vector field $V_2$ of rank $n$ has
constant coefficients if and only if: $[V_2,V_2]=0$. If $n$ is odd, then a
$2$-vector field $V_2$ of rank $n-1$ has constant coefficients if and only
if the Pfaffian system $P=\{ w\in T_x^\ast M:i_wV_2=0\} $ is integrable
and the restriction of $[V_2,V_2]$ to each integral submanifold of $P$, vanishes.
\end{corollary}
\begin{proof}
If $n$ is odd, then $\operatorname*{rank}P=1$, as $\operatorname*{rank}
V_2=n-1$. Hence, there exists a local section of $P$ with constant
coefficients, $\omega _1=\lambda _idu^i\in \Gamma(U,P)$, and any other
section of $P$ is of the form $f\omega_1$.

If $n$ is even, the mapping $V_2^\natural \colon T^\ast M\to TM $,
$V_2^{\natural}(w)=i_{w}V_2$, $\forall w\in T^\ast M$, is an isomorphism,
as $\operatorname*{rank}V_2=n$, and a $2$-form can be defined by the formula
$\omega_2(X,Y)=V_2((V_2^{\natural})^{-1}X,(V_2^{\natural})^{-1}Y)$,
$\forall X,Y\in TM$. According to \cite{Lichnerowicz} (also see
\cite[III.8.12--2]{Libermann-Marle}), $\omega_2$ is closed if and only if
$[V_2,V_2]=0$. In this case, by virtue of Darboux's theorem, there exist
coordinates $(p_i,q^i)$, $1\leq i\leq\frac{n}{2}$, such that
$\omega _2=dp_i\wedge dq^i$; hence
$V_2=\frac{\partial }{\partial p_i}\wedge\frac{\partial }{\partial q^i}$
has constant coefficients.

If $n$ is odd and $P$ is integrable, then the distribution $\ker
P\subset\mathfrak{X}(M)$ is involutive, since if $\omega_1$ is an arbitrary
local section of $P$, then $d\omega_1=\omega_1\wedge\xi_1$ for some
$1$-form $\xi_1$. Therefore, for all $X,Y\in \Gamma(U,\ker P)$ we have
\[
d\omega_1(X,Y)=\left( \omega_1\wedge\xi_1\right) (X,Y)=\omega
_1(X)\xi_1(Y)-\omega_1(Y)\xi_1(X)=0,
\]
since $\omega_1(X)=\omega_1(Y)=0$. Furthermore
\begin{align*}
0  & =d\omega_1(X,Y)=X(\omega_1(Y))-Y(\omega_1(X))-\omega_1([X,Y])\\
& =-\omega_1([X,Y]).
\end{align*}
Accordingly $[X,Y]\in \Gamma(U,\ker P)$. From Frobenius theorem, there exist
coordinates $x^1,\dotsc,x^n$ such that $\ker P=\left\langle \frac
{\partial }{\partial x^1},\dotsc,\frac{\partial }{\partial x^{n-1}}
\right\rangle $, so that in this system we have $P=\left\langle
dx^n\right\rangle $; hence $V_2=\sum_{1\leq i<j\leq n-1}F_{ij}
\frac{\partial }{\partial x^i}\wedge\frac{\partial }{\partial x^i}$, and one
concludes by applying the proof of the even-dimensional case to the
restriction of $V_2$ to each integral submanifold $x^n=\lambda$ of $P$.
\end{proof}
\begin{proposition}
\label{propos0.5}
If $V_{n-1}$ is a $(n-1)$-vector field such that
$(V_{n-1})_{x_0}\neq 0$ at a point $x_0\in M$, then
$V_{n-1}$ has constant coefficients on an open neigbourhood
$U$ of $x_0$ if and only if the following two conditions hold:
\begin{enumerate}
\item[\emph{(1)}] The Pfaffian system $P$ associated to $V_{n-1}$
according to Proposition \ref{propos0} satisfies
the properties therein stated.
\item[\emph{(2)}] Let $\eta $ be the $\wedge ^nTM$-valued differential
$1$-form defined by
\[
\begin{array}
[c]{llll}
\eta (X)= & X\wedge V_{n-1}, &  & \forall X\in \mathfrak{X}(M).
\end{array}
\]
There exists a derivation law
\[
\begin{array}
[c]{ll}
\nabla \colon \mathfrak{X}(M)\times\wedge ^n\mathfrak{X}(M)
\to \wedge ^n\mathfrak{X}(M), & (X,W_n)\mapsto \nabla _XW_n,
\end{array}
\]
with vanishing curvature (i.e., $R_\nabla =0$) such that
$d^\nabla \eta =0$.
\end{enumerate}
\end{proposition}
\begin{proof}
If $V_{n-1}$ has constant coefficients, i.e.,
$V_{n-1}=c_i\frac{\partial }{\partial x^1}\wedge \cdots \wedge
\widehat{\frac{\partial }{\partial x^i}}\wedge \cdots \wedge
\frac{\partial }{\partial x^n}$, then
$\eta=(-1)^{i-1}c_{i}dx^i\otimes
\frac{\partial }{\partial x^1}\wedge\cdots \wedge
\frac{\partial }{\partial x^n}$ and it suffices to consider
the derivation law given by
\[
\begin{array}
[c]{lll}
\nabla _X\left( f\frac{\partial}{\partial x^1}\wedge \cdots \wedge
\frac{\partial }{\partial x^n}\right)
=X(f)\frac{\partial }{\partial x^1}\wedge \cdots \wedge
\frac{\partial }{\partial x^n},
& \forall f\in C^\infty (M),
& \forall X\in \mathfrak{X}(M).
\end{array}
\]
Conversely, assume that items (1) and (2) hold.

Let $U$ be the neigbourhood of $x_0$ defined by
$\{ x\in M:(V_{n-1})_x\neq 0\} $. Shrinking $U$
if necessary we can assume that $U$ is a
coordinate domain. If $V_n^0$ is a basis for
$\Gamma (U,\wedge ^nTM)$, then there exists a unique
$\omega _1$ in $\Omega ^1(U)$ such that
$i_{\omega _1}V_n^0=V_{n-1}$; hence $i_{\omega _1}V_{n-1}=0$,
and thus $\operatorname*{rank}P_{|U}\geq 1 $. According to
the assumption in item (i) we conclude that there exists
a function $f$ in $C^{\infty}(U)$ such that
$(df)_x\neq 0$, $\forall x\in U$, and $i_{df}V_{n-1}=0$.
Again by shrinking $U$, we can assume that there exists
a coordinate system $(U;x^1,\dotsc,x^n)$ such that
$f=x^n$. We have
$V_{n-1}=F\frac{\partial }{\partial x^1}\wedge \cdots
\wedge\frac{\partial }{\partial x^{n-1}}$,
$\eta =\alpha \otimes V_n$, with $F\in C^\infty (U)$,
$\alpha =(-1)^{n-1}Fdx^n$ and
$V_n=\frac{\partial }{\partial x^1}\wedge \cdots
\wedge \frac{\partial }{\partial x^n}$.
The connection form $\omega $\ of $\nabla $ is defined by the formula
$\nabla _XV_n=\omega (X)V_n$, and the curvature is then given by
$R_\nabla =d\omega \otimes V_n$. As $R_\nabla =0$, we have
$\omega =d\rho $, and substituting $\exp (-\rho) V_n$ for $V_n$,
it follows $\nabla _X(V_n)=0$.
Therefore $d^\nabla \eta =d\alpha\otimes V_n
=(-1)^{n-1}dF\wedge dx^n \otimes V_n=0$, and, consequently, $F=F(x^n)$.
As $F(x_0)\neq 0$, the following formulas $y^1=\frac{x^1}{F}$,
$y^i=x^i$, $2\leq i\leq n$, determine a change of coordinates
and we have $V_{n-1}=\frac{\partial }{\partial y^1}\wedge
\cdots \wedge \frac{\partial }{\partial y^{n-1}}$.
\end{proof}
\begin{remark}
\label{remark3}
If $\nabla $ is a derivation law as in Proposition \ref{propos0.5}
and $D_\nabla $ is the differential operator generating
Schouten-Nijenhuis bracket introduced in \cite[\S 2]{Koszul},
then $(D_\nabla )^2=0$, as follows from the formula \cite[(2.2)]{Koszul}.
\end{remark}
\begin{remark}
\label{remark4}
For $q\geq \frac{1}{2}(n+2)$ the condition $[V_q,V_q]=0$
holds automatically, as $\deg [V_q,V_q]\geq 2q-1\geq n+1$.
\end{remark}
\begin{example}
\label{example6}
An $n$-vector field $V_n$ has constant coefficients on a
neighbourhood of $x\in M$ if and only if, either $V_n$ vanishes on a
neighbourhood of $x$, or $V_{n}$ does not vanish at $x$. In fact, in the
second case, there exists a unique $n$-form $\omega _n$ such that
 $\omega _n(V_n)=1$, and we conclude by virtue of Example \ref{example2}.
\end{example}
\subsection{Conformal constant coefficients}
Each volume form $\omega _n$ induces an isomorphism
\[
\begin{array}
[c]{l}
\iota _{pq}\colon \wedge ^qT(M)\to \wedge ^pT^\ast (M),
\quad n=p+q,\smallskip \\
\iota _{pq}\left( X_1\wedge \cdots \wedge X_q\right)
=i_{X_1\wedge \cdots \wedge X_q}\omega _n=i_{X_1}
\left( \ldots i_{X_q}\left( \omega _n\right) \right) ,\\
\forall X_1,\dotsc,X_q\in T_xM.
\end{array}
\]
If $\omega _n=\rho dx^1\wedge \cdots \wedge dx^n$, then for every
$I=(i_1,\dotsc,i_q)\in \mathcal{I}_q^n$, we have
\[
\begin{array}
[c]{rl}
\iota _{pq}
\left( \tfrac{\partial }{\partial x^{i_1}}\wedge \cdots
\wedge \tfrac{\partial }{\partial x^{i_q}}
\right)
= & i_{\frac{\partial }{\partial x^{i_1}}}
\Bigl(
\ldots i_{\frac{\partial }{\partial x^{i_q}}}\left( \omega _n\right)
\Bigr)
\\
= & \rho i_{\frac{\partial }{\partial x^{i_1}}}
\Bigl(
\ldots (-1)^{i_q-1}dx^1\wedge \cdots \wedge \widehat{dx^{i_q}}
\wedge \cdots \wedge dx^n
\Bigr)
\\
= & (-1)^{|I|-q}\rho dx^1\wedge \cdots \wedge \widehat{dx^{i_1}}
\wedge \cdots \wedge \widehat{dx^{i_q}}\wedge \cdots \wedge dx^n.
\end{array}
\]
Hence
\begin{multline*}
\iota _{pq}\left(
{\textstyle\sum\nolimits_{I\in \mathcal{I}_q^n}}
F_I\tfrac{\partial }{\partial x^{i_1}}\wedge \cdots
\wedge \tfrac{\partial }{\partial x^{i_q}}\right) \\
=\rho
{\textstyle\sum\nolimits_{I\in \mathcal{I}_q^n}}
(-1)^{|I|-q}F_Idx^1\wedge \cdots \wedge\widehat{dx^{i_1}}\wedge
\cdots \wedge \widehat{dx^{i_q}}\wedge \cdots \wedge dx^n.
\end{multline*}
Dually, an isomorphism can also be defined
\[
\begin{array}
[c]{l}
\iota _{qp}^\ast \colon \wedge ^pT^\ast (M)\to \wedge ^qT(M),
\quad
n=p+q,\smallskip \\
\iota _{qp}^\ast \left( w_1\wedge \cdots \wedge w_p\right)
=i_{w_1\wedge \cdots \wedge w_p}V_n
=i_{w_1}\left( \ldots i_{w_p}\left( V_n\right) \right) ,\\
\forall w_1,\dotsc,w_p\in T_x^\ast M,
\end{array}
\]
$V_n$ being an $n$-vector field that does not vanish at any point.
\begin{definition}
\label{def2}
A $q$-vector field $V_q$ (resp.\ a differential $p$-form
$\omega _p$) is said to have \emph{conformal constant coefficients} if there
exists a $q$-vector field $V_q^\prime $ (resp.\ a differential $p$-form
$\omega _p^\prime $) with constant coefficients and function
$f\in  C^\infty (M)$ such that $V_q=fV_q^\prime $ (resp.\
$\omega _p=f\omega _p^\prime $).
\end{definition}
\begin{remark}
\label{remark5}
From the previous formulas it follows that if $V_q$ has
conformal constant coefficients, then the $p$-form $\iota _{pq}(V_q)$ has
also conformal constant coefficients. In fact, if $V_q^\prime $ has
constant coefficients on $(u^1,\dotsc,u^n)$, then $\omega _n=\rho
du^1\wedge \cdots \wedge du^n$; hence $\iota _{pq}(V_q)=f\rho
i_{V_q^\prime }\left( du^1\wedge \cdots \wedge du^n\right) $, which is
a $p$-form with conformal constant coefficients.
\end{remark}
\begin{remark}
\label{remark6}
If $V_q$ and $\omega _n$ have constant coefficients in the
same coordinate system, then $i_{V_q}\omega _n$ has also constant
coefficients in that system. If $n=pk$, $\omega _p$ has constant coefficients
and $\omega _n=\omega _p\wedge\overset{(k}{\ldots}\wedge\omega _p$ is a
volume form, then the $q$-vector field $V_q$, $q=p(k-1)$, defined by
$i_{V_q}\omega _n=\omega _p$, has also constant coefficients.
\end{remark}
\begin{example}
\label{example7}
An $i$-form $\omega _i$ on $M$ with $i=1$ or $i=2$ of constant class
and such that $(\omega _i)_x\neq 0$ at point $x\in M$ has conformal
constant coefficients if and only if $\omega _i\wedge d\omega _i=0$.
In fact, if $\omega _ i=\rho \omega^\prime _i$, with
$\rho \in C^\infty (M)$ and $\omega ^\prime _i$ has constant
coefficients, then $d\omega _ i=d\rho \wedge \omega^\prime _i$,
and hence $\omega _ i\wedge d\omega _ i=0$. Conversely, if
the previous equation holds, then the class of $\omega _i$
is $2$ and by virtue of Darboux's theorem
(cf.\ \cite[VI, \S 4.1]{Godbillon}) there exists a coordinate
system $(x^i)_{i=1}^n$ such that either $\omega _1=(1+x^1)dx^2$
or $\omega _2=\rho dx^1\wedge dx^2$.
\end{example}
\begin{example}
\label{example7.5}
If $V_{n-1}$ is a $(n-1)$-vector field on $M$ such that
$(V_{n-1})_{x_0}\neq 0$ on a point $x_{0}\in M$
and the Pfaffian system $P$ associated to $V_{n-1}$
according to Proposition \ref{propos0} satisfies the properties
therein stated, then $V_{n-1}$ has conformal constant coefficients.
This is proved by proceeding as in the proof of Proposition \ref{propos0.5}.
\end{example}
The results of this section prove that for very high or very low
 degrees $p$ or $q$, differential $p$-forms
and $q$-vector fields with constant coefficients admit
simple geometric characterizations obtained by using
classical theorems.

For intermediate degrees the situation is different
and one must resort to consider partial differential
systems of geometric nature that allow to characterize
differential forms and vector fields with constant
coefficients. This is the purpose of the next two sections.
\section{The associated linear connection}
Given a local basis $\mathbf{f}=(X_1,\dotsc,X_n)$ of the bundle of linear
frames $F(M)$, or equivalently, $(X_1,\dotsc,X_n)$
is a local basis of the module $\mathfrak{X}(M)$,
let $\nabla ^{\mathbf{f}}$ denote the only linear
connection on $M$ parallelizing $\mathbf{f}$, i.e.,
$\nabla _{X_i}^{\mathbf{f}}X_j=0$, $\forall i,j=1,\dotsc,n$.
If $\mathbf{u}=(u^1,\dotsc,u^n)$ is a coordinate system on $M$, we write
$\nabla ^{\mathbf{u}}=\nabla ^{\mathbf{f}}$, with
$X_i=\frac{\partial }{\partial u^i}$,
$1\leq i\leq n$. The curvature of $\nabla ^{\mathbf{f}}$ vanishes (i.e.,
$\nabla ^{\mathbf{f}}$ is flat) and its torsion vanishes if and only if there
exist coordinates such that $\nabla ^{\mathbf{f}}=\nabla ^{\mathbf{u}}$.
\begin{proposition}
\label{propos1}
A $p$-form $\omega _p$ (resp.\  a $q$-vector field $V_q$) has
constant coefficients on the coordinate system $\mathbf{u}=(u^1,\dotsc,u^n)$
if and only if: $\nabla ^{\mathbf{u}}\omega _p=0$ (resp.\
$\nabla ^{\mathbf{u}}V_q=0$).

In other words, $\omega _p$ (resp.\ $V_q$) has constant coefficients
if and only if there exists a symmetric linear connection $\nabla $ on $M$
such that
$R^\nabla =0$, $\nabla \omega _p=0$ (resp.\ $\nabla V_q=0$).
\end{proposition}
\begin{proof}
The condition $\nabla ^{\mathbf{u}}\omega _p=0$ means
$\nabla _X^{\mathbf{u}}\omega _p=0$, $\forall X\in \mathfrak{X}(M)$.
We have $\nabla ^{\mathbf{u}}(du^J)=0$, as
\[
\Bigl(
\nabla _{\tfrac{\partial }{\partial u^i}}^{\mathbf{u}}(du^J)
\Bigr)
\tfrac{\partial }{\partial u^k}=\tfrac{\partial }{\partial u^i}
\left( du^J\left( \tfrac{\partial }{\partial u^k}\right) \right)
-du^J
\Bigl(
\nabla _{\frac{\partial }
{\partial u^i}}^{\mathbf{u}}\tfrac{\partial }{\partial u^k}
\Bigr)
=0.
\]
Consequently, if $\omega _p=\sum_{I\in \mathcal{I}_p^n}F_Idu^{I}$ on the
system $\mathbf{u}=(u^1,\dotsc,u^n)$, then the equation $\nabla
_X^{\mathbf{u}}\omega _p=0$ is equivalent to saying $X(F_I)=0$, $\forall
X\in \mathfrak{X}(M)$, and all multi-indices $I\in \mathcal{I}_p^n$; in
other words, all the coefficients of $\omega _p$ are constant.

The proof for the contravariant case is analogous.
\end{proof}
\begin{lemma}
\label{lemma1}
If $\nabla _{\frac{\partial }{\partial x^i}}^{\mathbf{u}}
\Bigl(
\frac{\partial }{\partial x^j}
\Bigr)
=\Gamma _{ij}^h\frac{\partial }{\partial x^h}$, then
$\Gamma _{ij}^b=v_h^b\frac{\partial^2u^h}{\partial x^i\partial x^j}$,
where $(v_j^i)$ denotes the inverse matrix of the Jacobian matrix
$\left( \frac{\partial u^i}{\partial x^j}\right) $, $b,h,i,j=1,\dotsc,n$.
\end{lemma}
\begin{proof}
From $\frac{\partial }{\partial x^j}
=\frac{\partial u^h}{\partial x^j}\frac{\partial }{\partial u^h}$,
it follows:
\[
\nabla _{\frac{\partial }{\partial x^i}}^{\mathbf{u}}
\left( \tfrac{\partial }{\partial x^j}\right)
=\tfrac{\partial^2u^h}{\partial x^i\partial x^j}
\tfrac{\partial }{\partial u^h}
+\tfrac{\partial u^h}{\partial x^j}
\nabla _{\frac{\partial }{\partial x^i}}^{\mathbf{u}}
\Bigl(
\tfrac{\partial }{\partial u^h}
\Bigr) .
\]
As the second term in the right-hand side vanishes, taking the formulas
$\tfrac{\partial u^a}{\partial x^b}v_c^b=\delta _c^a$,
$v_c^b\tfrac{\partial }{\partial x^b}
=\tfrac{\partial u^a}{\partial x^b}v_c^b
\tfrac{\partial }{\partial u^a}
=\delta _c^a
\tfrac{\partial }{\partial u^a}
=\tfrac{\partial }{\partial u^c}$ into account, we have
\[
\nabla _{\frac{\partial }{\partial x^i}}^{\mathbf{u}}
\left( \tfrac{\partial }{\partial x^j}\right)
=\tfrac{\partial ^2u^h}{\partial x^i\partial x^j}
\tfrac{\partial }{\partial u^h}=v_h^b
\tfrac{\partial ^2u^h}{\partial x^i\partial x^j}
\tfrac{\partial }{\partial x^b},
\]
thus concluding.
\end{proof}
\section{The associated partial differential system}
In Proposition \ref{propos1} we have proved that
a $p$-form or a $q$-vector field have constant
coefficients if and only if there exists a linear
connection of vanishing curvature parallelizing them.
Below, this condition is proved to be equivalent
to the existence of solution to a partial differential
system of equations, which we describe explicitly.
\begin{theorem}
\label{th1}
Let $n=\dim M$ and let $O_{M}^2\subset J^2(M,\mathbb{R}^n)$
be the open subbundle of $2$-jets $j_x^2(u^1,\dotsc,u^n)$
such that
$(du^1\wedge \cdots \wedge du^n)_x\neq 0$.

Given a differential $p$-form and $q$-vector field,
\[
\omega _p=\sum_{I\in \mathcal{I}_p^n}F_Idx^{I},\quad V_q
=\sum _{I\in \mathcal{I}_q^n}\bar{F}_I\partial x^{I},
\]
on $M$, let $\Phi _{j,J}$ (resp.\  $\bar{\Phi }_{j,J}$)
be the functions defined by the following formulas:
\begin{equation}
\begin{array}
[c]{ll}
\sum_{J\in \mathcal{I}_p^n}\Phi _{j,J}dx^J
= & \sum\limits_{I\in \mathcal{I}_p^n}\sum\limits_{h=1}^p
\sum\limits_{i=1}^n\Gamma _{ji}^{i_h}F_Idx^{i_1}\wedge
\cdots\wedge
\underset{(i_h}{\underline{dx^i}}\wedge \cdots \wedge dx^{i_p},\\
\sum_{J\in \mathcal{I}_q^n}\bar{\Phi}_{j,J}\partial x^J
= & \sum\limits_{I\in \mathcal{I}_q^n}\sum\limits_{h=1}^{q}
\sum\limits_{i=1}^n\Gamma _{ji}^{i_h}
\bar{F}_I\partial x^{i_1}\wedge \cdots \wedge
\underset{(i_h}{\underline{\partial x^i}}\wedge \cdots
\wedge\partial x^{i_q},\\
1\leq j\leq n, &
\end{array}
\label{derivada1}
\end{equation}
where $\Gamma _{bc}^a$ are the functions introduced in
\emph{Lemma \ref{lemma1}} and the notations
$\underset{(i_h}{\underline{dx^i}}$ and
$\underset{(i_h}{\underline{\partial x^i}}$ mean
that the $1$-form $dx^i$ and the vector field
$\partial x^i$ must be inserted in the $i_h$-th place
of the exterior product.

The form $\omega _p$ (resp.\  the vector field $V_q$) has constant
coefficients on a neighbourhood of $x\in M$ if and only if the following
second-order partial differential system defined on $J^2(M,\mathbb{R}^n)$:
\begin{equation}
\left.
\begin{array}
[c]{rll}
\text{\emph{(i)}} & \tfrac{\partial F_J}{\partial x^j}=\Phi _{j,J},
& \forall J\in \mathcal{I}_p^n\medskip \\
\text{\emph{(ii)}}
& \tfrac{\partial\bar{F}_J}{\partial x^j}=-\bar{\Phi }_{j,J},
& \forall J\in \mathcal{I}_q^n
\end{array}
\right] \;1\leq j\leq n,
\label{ecs1}
\end{equation}
\noindent\medskip admits a local solution $(u^1,\dotsc,u^n)$ defined
on a neighbourhood of $x$ in $O_M^2$.
\end{theorem}
\begin{proof}
According to Proposition \ref{propos1} we must impose
$\nabla ^{\mathbf{u}}\omega _p=0$, i.e.,
\[
\nabla _{\frac{\partial }{\partial x^j}}^{\mathbf{u}}\omega _p
=\sum \limits_{I\in \mathcal{I}_p^n}
\Bigl(
\tfrac{\partial F_I}{\partial x^j}dx^{I}
+\sum\limits_{h=1}^pF_Idx^{i_1}\wedge \cdots \wedge
\nabla _{\frac{\partial }{\partial x^j}}^{\mathbf{u}}
\left( dx^{i_h}\right) \wedge \cdots \wedge dx^{i_p}
\Bigr) ,
\]
where $I=(i_1,\dotsc,i_p)$. Moreover, we have
\[
\nabla _{\frac{\partial }{\partial x^j}}^{\mathbf{u}}
\left( dx^{i_h}\right) \left( \tfrac{\partial }{\partial x^i}\right)
=\tfrac{\partial }{\partial x^j}
\left( dx^{i_h}\left( \tfrac{\partial }{\partial x^i}\right) \right)
-dx^{i_h}\left(
\nabla _{\frac{\partial }{\partial x^j}}^{\mathbf{u}}
\tfrac{\partial }{\partial x^i}
\right)
=-\Gamma _{ji}^{i_h}.
\]
Hence $\nabla _{\frac{\partial }{\partial x^j}}^{\mathbf{u}}
\left( dx^{i_h}\right) =-\Gamma _{ji}^{i_h}dx^i$, and thus
\begin{align*}
0  & =\nabla _{\frac{\partial }{\partial x^j}}^{\mathbf{u}}\omega _p\\
& =\sum\nolimits_{I\in \mathcal{I}_p^n}
\left(
\tfrac{\partial F_I}{\partial x^j}dx^I
-\sum\nolimits_{h=1}^p\Gamma _{ji}^{i_h}
F_Idx^{i_1}\wedge \cdots \wedge
\underset{(i_h}{\underline{dx^i}}\wedge \cdots \wedge dx^{i_p}\right) ,
\end{align*}
which allows one to conclude.

The proof of the contravariant case is analogous, replacing the formulas
$\nabla _{\frac{\partial }{\partial x^j}}^{\mathbf{u}}\left( dx^h\right)
=-\Gamma _{ji}^hdx^i$ by
$\nabla _{\frac{\partial }{\partial x^j}}^{\mathbf{u}}
\left( \frac{\partial }{\partial x^h}\right)
=\Gamma _{jh}^i\frac{\partial }{\partial x^i}$.
\end{proof}
The right-hand side of \eqref{derivada1} can be computed
in order to obtain the function $\Phi _{j,J}$ explicitely,
by simply applying the rules of calculation in the exterior algebra.
\begin{example}
As stated in Remark 1.9, the first case not included in the examples in
section 1.1 is obtained for the values $p=3$, $n=5$. In this case, we have
\[
\begin{array}
[c]{rl}
\omega_3= & F_{123}dx^1\wedge dx^2\wedge dx^3+F_{124}dx^1\wedge
dx^2\wedge dx^4+F_{125}dx^1\wedge dx^2\wedge dx^5\\
& +F_{134}dx^1\wedge dx^3\wedge dx^4+F_{135}dx^1\wedge dx^3\wedge
dx^5+F_{145}dx^1\wedge dx^4\wedge dx^5\\
& +F_{234}dx^2\wedge dx^3\wedge dx^4+F_{235}dx^2\wedge dx^3\wedge
dx^5+F_{245}dx^2\wedge dx^4\wedge dx^5\\
& +F_{345}dx^3\wedge dx^4\wedge dx^5.\\
&
\end{array}
\]
For simplicity's sake we write $F_{hij}=0$ if there is a repeated index in
$(h,i,j)$; otherwise $F_{hij}=\varepsilon F_{abc}$, with $a<b<c$ and
$\varepsilon$ is the sign of the permutation $a\mapsto h$, $b\mapsto i$,
$c\mapsto j$. Expanding on the first formula of item (1) in Theorem 3.1, we
have
\[
\begin{array}
[c]{l}
\Phi_{j,123}dx^1\wedge dx^2\wedge dx^3+\Phi_{j,124}dx^1\wedge
dx^2\wedge dx^4+\Phi_{j,125}dx^1\wedge dx^2\wedge dx^5\\
+\Phi_{j,134}dx^1\wedge dx^3\wedge dx^4+\Phi_{j,135}dx^1\wedge
dx^3\wedge dx^5+\Phi_{j,145}dx^1\wedge dx^4\wedge dx^5\\
+\Phi_{j,234}dx^2\wedge dx^3\wedge dx^4+\Phi_{j,235}dx^2\wedge
dx^3\wedge dx^5+\Phi_{j,245}dx^2\wedge dx^4\wedge dx^5\\
+\Phi_{j,345}dx^3\wedge dx^4\wedge dx^5=\medskip\\
\sum\limits_{I\in\mathcal{I}_3^5}\sum\limits_{i=1}^5\Gamma_{ji}^{i_1
}F_{i_1i_2i_3}dx^{i}\wedge dx^{i_2}\wedge dx^{i_3}+\sum
\limits_{I\in\mathcal{I}_3^5}\sum\limits_{i=1}^5\Gamma_{ji}^{i_2
}F_{i_1i_2i_3}dx^{i_1}\wedge dx^{i}\wedge dx^{i_3}\\
\multicolumn{1}{r}{+\sum\limits_{I\in\mathcal{I}_3^5}\sum\limits_{i=1}^5
\Gamma_{ji}^{i_3}F_{i_1i_2i_3}dx^{i_1}\wedge dx^{i_2}\wedge
dx^i.}
\end{array}
\]
Hence
\[
\begin{array}
[c]{rl}
\Phi_{j,123}= & \Gamma_{j1}^h\left(  F_{h23}-F_{2h3}+F_{23h}\right)
+\Gamma_{j2}^h\left( F_{1h3}-F_{h13}-F_{13h}\right) \\
& \multicolumn{1}{r}{+\Gamma_{j3}^h\left(  F_{h12}-F_{1h2}+F_{12h}\right) ,}\\
\Phi_{j,124}= & \Gamma_{j1}^h\left(  F_{h24}-F_{h24}+F_{24h}\right)
+\Gamma_{j2}^h\left( F_{1h4}-F_{h14}-F_{14h}\right) \\
& \multicolumn{1}{r}{+\Gamma_{j4}^h\left(  F_{h12}-F_{1h2}+F_{12h}\right) ,}\\
\Phi_{j,125}= & \Gamma_{j1}^h\left( F_{h25}-F_{h25}+F_{25h}\right)
+\Gamma_{j2}^h\left( F_{1h5}-F_{h15}-F_{15h}\right) \\
& \multicolumn{1}{r}{+\Gamma_{j5}^h\left( F_{h12}-F_{1h2}+F_{12h}\right) ,}\\
\multicolumn{1}{c}{\vdots} & \multicolumn{1}{c}{\vdots}
\end{array}
\]
and thus analogously for the remaining components.

Moreover, the cases studied in section 1.2 (the contravariant case) are:
$q=1$ (Example 1.13), $q=2$ and $V_2$ being of maximum rank (Corollary
1.15), $q=n-1$ (Proposition 1.16), and $q=n$ (Example 1.19) with arbitrary
$n$. Therefore, the first case not included in the examples in section 1.1 is
also obtained for the values $q=3$, $n=5$. Consequently, replacing
$\Gamma_{ij}^h$ by $-\Gamma_{jh}^i$ in the previous formulas, the
corresponding values for $\bar{\Phi }_{j,abc}$, $1\leq a<b<c\leq5$,
are also obtained.
\end{example}
\begin{corollary}
\label{cor1}
Let $\Psi_{j,J}$, $J\in \mathcal{I}_p^n$ (resp.\
$\bar{\Psi }_{j,J}$,
$J\in \mathcal{I}_q^n$) be the function obtained by replacing
$\Gamma _{cd}^a=v_b^a\frac{\partial ^2u^b}{\partial x^c\partial x^d}$
by $\Gamma _{cd}^a=-\tfrac{\partial v_b^a}{\partial x^d}
\left( V^{-1}\right) _c^b$ in the formulas above, where
$V=(v_j^i)$ is the inverse matrix
of the Jacobian matrix $\left( \frac{\partial u^i}{\partial x^j}\right) $.
The systems \emph{\eqref{ecs1}-(i)} and\emph{\eqref{ecs1}-(ii)} of PDS of
second order of \emph{Theorem \ref{th1}}, defined in $O_M^2$, are
respectively equivalent to the following systems of first order in the
variables $v_j^i$, $i,j=1,\dotsc,n$, both defined on
$J^1(M,GL(n,\mathbb{R}))$:
\begin{equation}
\begin{array}
[c]{lll}
\tfrac{\partial F_J}{\partial x^j}=\Psi_{j,J},
& \forall J\in \mathcal{I}_p^n,
& 1\leq j\leq n,\smallskip \\
v_{m}^j\tfrac{\partial v_l^k}{\partial x^j}
=\tfrac{\partial v_m^k}{\partial x^i}v_l^i
& 1\leq k\leq n,1\leq l<m\leq n,
&
\end{array}
\label{ecs2}
\end{equation}
\begin{equation}
\begin{array}
[c]{lll}
\tfrac{\partial\bar{F}_J}{\partial x^j}=\bar{\Psi }_{j,J},
& \forall J\in \mathcal{I}_q^n, & 1\leq j\leq n,
\smallskip \\
v_m^j\tfrac{\partial v_l^k}{\partial x^j}
=\tfrac{\partial v_{m}^k}{\partial x^i}v_l^i
& 1\leq k\leq n,1\leq l<m\leq n. &
\end{array}
\label{bar_ecs2}
\end{equation}
\end{corollary}
\begin{proof}
By taking derivatives on
$v_b^a\frac{\partial u^b}{\partial x^c}=\delta _c^a$
with respect to $x^d$ it follows:
$\tfrac{\partial v_b^a}{\partial x^d}
\tfrac{\partial u^b}{\partial x^c}
+v_b^a\tfrac{\partial ^2u^b}{\partial x^c\partial x^d}=0$.
According to Lemma \ref{lemma1}:
$\Gamma _{cd}^a=v_b^a\frac{\partial ^2u^b}{\partial x^c\partial x^d}$;
hence
\[
\Gamma _{cd}^a=-\tfrac{\partial v_b^a}{\partial x^d}
\tfrac{\partial u^b}{\partial x^c}
=-\tfrac{\partial v_b^a}{\partial x^d}
\left( V^{-1}\right) _c^b.
\]
The system \eqref{ecs1} of Theorem \ref{th1} is thus of first order with
respect to the variables $v_j^i$, and hence, it is defined on
$J^1(M,GL(n,\mathbb{R}))$; but this is not enough to solve it by taking the
functions $v_j^i$ as unknowns, since it is not ensured that the inverse
matrix $V^{-1}$ is a Jacobian matrix. For this, the following equations must
also be verified: $\tfrac{\partial(V^{-1})_i^h}{\partial x^j}
=\tfrac{\partial(V^{-1})_j^h}{\partial x^i}$, $h,i,j=1,\dotsc,n$,
or equivalently,
\begin{equation}
(V^{-1})_a^h\tfrac{\partial v_b^a}{\partial x^j}(V^{-1})_i^b
=(V^{-1})_a^h\tfrac{\partial v_b^a}{\partial x^i}(V^{-1})_j^b,
\quad h,i,j=1,\dotsc,n,
\label{derivada2}
\end{equation}
as follows by taking derivatives in the equality $v_b^a(V^{-1})_c^b
=\delta _c^a$ with respect to $x^d$ and solving it for
$\frac{\partial(V^{-1})_c^b}{\partial x^d}$. Finally, by multiplying
the two sides of \eqref{derivada2} by $v_h^k$ and by summing up
over $h$, it follows: $\tfrac{\partial v_b^k}{\partial x^j}(V^{-1})_i^b
=\tfrac{\partial v_b^k}{\partial x^i}(V^{-1})_j^b$, and by
multiplying both sides by $v_m^i$ and summing up over $i$, we thus obtain
\[
\tfrac{\partial v_l^k}{\partial x^j}
=\tfrac{\partial v_b^k}{\partial x^i}(V^{-1})_j^bv_l^i.
\]
Furthermore, by multiplying both sides by $v_m^j$ and summing up over $j$,
we conclude.
\end{proof}
\begin{corollary}
\label{cor2}
The integrability conditions of the system \emph{\eqref{ecs1}-(i)}
(resp.\ \emph{\eqref{ecs1}-(ii)}) with respect to the unknowns $F_J$ (resp.\
$\bar{F}_J$) hold identically.

Accordingly, if $M$ is of class $C^\omega $, given a point $x_{0}\in M$,
scalars $\lambda _J$ (resp.\  $\bar{\lambda }_J$), $1\leq j\leq n$, and a
coordinate system $(u^1,\dotsc,u^n)$ of class $C^\omega $ on $M$, there
exists a unique differential $p$-form $\omega _p$\ (resp.\ a $q$-vector field
$V_q$) on $M$ with constant coefficients in the system $(u^1,\dotsc,u^n)$
such that $F_J(x_0)=\lambda _J$ (resp.\  $\bar{F}_J(x_0)=\bar{\lambda }_J$).
\end{corollary}
\begin{proof}
From the equations \eqref{ecs1}-(i) it follows:
\[
\begin{array}
[c]{ll}
\tfrac{\partial ^2F_J}{\partial x^j\partial x^l}
=\tfrac{\partial \Phi _{j,J}}{\partial x^l},
& \tfrac{\partial ^2F_J}{\partial x^l\partial x^j}
=\tfrac{\partial \Phi _{l,J}}{\partial x^j}.
\end{array}
\]
Hence, integrability conditions of the system \eqref{ecs1}-(i) with respect to
the functions $F_J$ are $0=\tfrac{\partial \Phi _{j,J}}{\partial x^l}
-\tfrac{\partial \Phi _{l,J}}{\partial x^j}$.

As $\nabla ^{\mathbf{u}}$ is flat, we have
$\frac{\partial \Gamma _{bd}^a}{\partial x^c}
-\frac{\partial\Gamma _{bc}^a}{\partial x^d}
=\Gamma _{bc}^e\Gamma _{de}^a-\Gamma _{bd}^e\Gamma _{ce}^a$,
and as a calculation shows, the integrability conditions written
above hold identically.

The proof in the contravariant case is analogous.
\end{proof}
\begin{remark}
\label{remark6.5}
Certainly a $p$-form $\omega _p$ or a $q$-vector field $V_q$
have constant coefficients in an open subset $U$ if there exists
a coordinate system $(u^i)_{i=1}^n$ defined on $U$ such that
$L_{\frac{\partial }{\partial u^i}}\omega _p=0$ or
$L_{\frac{\partial }{\partial u^i}}V _q=0$, respectively,
for $1\leq i\leq n$.
By directly imposing these conditions, systems of equations
equivalent to \eqref{ecs1} are obtained. However, we prefer
the torsion-free flat linear connection method because
it reveals the geometry underlying the problem considered
and is closer to the classical formulations given in Geometry;
for example, a pseudo-Riemannian metric has constant
coefficients if and only if the curvature of its associated
Levi-Civita connection vanishes.
\end{remark}
\section{Computational aspects}
\begin{remark}
\label{remark7}
For all $n\geq 7$, $3\leq p\leq n-3$, the system \eqref{ecs2} has
$n\tbinom{n}{p}+n\tbinom{n}{2}$ equations in the $n^2+n^3$ variables
$v_b^a$, $\tfrac{\partial v_b^a}{\partial x^d}$, and we have
$n\tbinom{n}{p}+n\tbinom{n}{2}\geq n^2+n^3$. In fact, only
for $3\leq p\leq 4$, $n=7$ the inequality above turns into
an equality, i.e., $7\tbinom{7}{p}+7\tbinom{7}{2}=7^2+7^3$.
Accordingly, from dimension $8$ there are more equations
than unknowns and the system \eqref{ecs2} is generically
compatible and determined, so that the connection
$\nabla ^{\mathbf{u}}$, if exists, is unique.
\end{remark}
\begin{remark}
\label{remark8}
The systems (i) and (ii) in \eqref{ecs1} can be viewed as two
linear systems in the $n^3$ unknowns $\Gamma _{bc}^a$ each with
$n\tbinom{n}{p}$ equations.

Moreover, as a computation shows, we have
\[
\begin{array}
[c]{ll}
n\tbinom{n}{p}<n^3, & \text{if }2\leq n\leq 7
\text{ and }2\leq p\leq n-2,\smallskip\\
n\tbinom{n}{p}<n^3, & \text{if }n=8
\text{ and }p=2,3,5,6,\smallskip\\
n\tbinom{n}{p}\geq n^3, & \text{if }
n=8\text{ and }p=4,\smallskip\\
n\tbinom{n}{p}<n^3, & \text{if }
n\geq 9\text{ and }n-2\leq p\leq n,
\smallskip\\
n\tbinom{n}{p}\geq n^3, & \text{if }
n\geq 9\text{ and }3\leq p\leq n-3.
\end{array}
\]
\end{remark}
\begin{remark}
\label{remark9}
Let $\mathcal{M}(x)$ be the coefficient matrix
of \eqref{ecs1}-(i)  at a point $x\in M$, and let
$\mathcal{M}^\prime (x)$ be its augmented matrix.

We have
$\operatorname*{rank}\mathcal{M}(x)
\leq \operatorname*{rank}\mathcal{M}^\prime (x)
\leq n^3$, and from Rouch\'e-Capelli theorem
it follows:
\begin{itemize}
\item
If $\operatorname*{rank}\mathcal{M}(x)
<\operatorname*{rank}\mathcal{M}^\prime (x)$,
then $\operatorname*{rank}\mathcal{M}(x^\prime )
<\operatorname*{rank}\mathcal{M}^\prime (x^\prime )$ for every
$x^\prime $ in a neighbourhood of $x$, and the $p$-form $\omega _p
=\sum_{I\in \mathcal{I}_p^n}F_Idx^I$ does not have constant
coefficients on a neighbourhood of $x$

A similar result also holds for a $q$-vector field.
\item
If $\operatorname*{rank}\mathcal{M}(x^\prime )=
\operatorname*{rank}\mathcal{M}^\prime (x^\prime )$
for every $x^\prime $ in a neighbourhood of $x$,
then the system \eqref{ecs1}-(i) admits at least a
solution $\Gamma _{bc}^a$, $a,b,c=1,\dotsc,n$;
but we need to impose that the curvature and torsion fields
of the linear connection $\nabla $ given by the formula
$\nabla _{\frac{\partial }{\partial x^b}}
\frac{\partial }{\partial x^c}=\Gamma _{bc}^a
\frac{\partial }{\partial x^a}$\ should vanish, or equivalently,
such that $\nabla =\nabla ^{\mathbf{u}}$ for some coordinate system
$\mathbf{u}=(u^1,\dotsc,u^n)$.

In other words, on the solution to \eqref{ecs1}-(i) we must impose
the following conditions:
\[
\begin{array}
[c]{rrl}
\text{(C)} & \frac{\partial \Gamma _{bd}^a}{\partial x^c}
-\frac{\partial \Gamma _{bc}^a}{\partial x^d}
= & \Gamma _{bc}^e\Gamma _{de}^a
-\Gamma _{bd}^e\Gamma _{ce}^a,\smallskip \\
\text{(T)} & \Gamma _{cb}^a= & \Gamma _{bc}^a,\medskip \\
& a,b,c,d,e= & 1,\dotsc,n.
\end{array}
\]
\end{itemize}
\end{remark}
\begin{remark}
\label{remark9.5}
If $n\binom{n}{p}\geq n^3$, $3\leq p\leq n-3$,
as in Remark \ref{remark8} and furthermore we have
$\operatorname*{rank}\mathcal{M}(x^\prime )=
\operatorname*{rank}\mathcal{M}^\prime (x^\prime )=n^3$,
then the solution to the system \eqref{ecs1}-(i)
can be written in the form $\Gamma _{bc}^a=Q_{bc}^a$,
$a,b,c=1,\dotsc,n$, where the functions $Q_{bc}^a$
are rational fractions of the functions $F_J$,
$\frac{\partial F_J}{\partial x^j}$,
$J\in \mathcal{I}_p^n$, $1\leq j\leq n$.
The contravariant case is completely similar. Hence,
in this case, there exist two second-order differential operators
\[
D\colon \mathcal{O}\subset \wedge ^p\Omega (\mathbb{R}^n)
\to C^\infty (\mathbb{R}^n)^k,
\quad
\bar{D}\colon \bar{\mathcal{O}}\subset
\wedge ^q\mathfrak{X}(\mathbb{R}^n)
\to C^\infty (\mathbb{R}^n)^{\bar{k}},
\]
for some integers $k$ and $\bar{k}$, where $\mathcal{O}$
and $\bar{\mathcal{O}}$ are dense open subsets, such that
a differential $p$-form $\omega _p\in \mathcal{O}$
and a $q$-vector field $V_q\in \bar{\mathcal{O}}$
have constant coefficients if and only if
$D(\omega _p)=0$ and $\bar{D}(V_q)=0$, respectively.

Although the operators $D$ and $\bar{D}$ can be computed feasibly,
since their computation is reduced to Linear Algebra operations,
their expressions become increasingly longer when the value of $n$
increases.

If $\operatorname*{rank}\mathcal{M}(x^\prime )=
\operatorname*{rank}\mathcal{M}^\prime (x^\prime )$
does not reach its maximum value, the order of the operators $D$
and $\bar{D}$ may be higher than $2$, if they exist.
\end{remark}
\begin{remark}
\label{remark10}
The case of a $(n-1)$-vector field
\[
V_{n-1}=\sum\nolimits_{i=1}^nF_i\tfrac{\partial}{\partial x^1}
\wedge \cdots \wedge \widehat{\tfrac{\partial }{\partial x^i}}
\wedge \cdots \wedge \tfrac{\partial }{\partial x^n}
\]
is exceptional because from Remark \ref{remark8} we know
that the system (ii) in \eqref{ecs1} is a linear system
in the $n^3$ unknowns $\Gamma _{bc}^a$ with $n\tbinom{n}{n-1}=n^2$
equations. In fact, taking the symmetry
$\Gamma _{bc}^a=\Gamma _{cb}^a$ into account, it follows that
the number of unknowns in this case is $\frac{1}{2}n^2(n+1)>n^2$.
Hence the system (ii)-\eqref{ecs1} is underdetermined and the
associated torsion-free flat linear connection is not unique.

With the previous notations, the system (ii)-\eqref{ecs1} can be
written as follows:
\begin{align}
\tfrac{\partial F_l}{\partial x^j}
& =\left( \delta _{l,1}-1\right)
\sum\nolimits_{h=1}^{l-1}
\Gamma _{jh}^hF_l+\left( \delta _{l,n}-1\right)
\sum\nolimits_{h=l+1}^n\Gamma _{jh}^hF_l
\label{E[j,l]}\\
& +\left( 1-\delta _{i,1}\right) \left( 1-\delta _{l,n}\right)
\sum\nolimits_{i=l+1}^n(-1)^{i-l}\Gamma _{jl}^iF_i
\nonumber\\
& +\left( 1-\delta _{i,n}\right) \left( 1-\delta_{l,1}\right)
\sum\nolimits_{i=1}^{l-1}(-1)^{l-i}\Gamma_{jl}^{i}F_i,
\nonumber\\
j,l & =1,\dotsc,n.
\nonumber
\end{align}
As a computation shows, on the dense open subset defined by
$F_1\cdots F_n\neq 0$ the rank of the $n^2\times \frac{1}{2}n^2(n+1)$
coefficient matrix of \eqref{E[j,l]} is $n^2$, i.e., the rank is maximum.
Hence $n^2$ of the unknowns $\Gamma _{bc}^a$, $a,b,c=1,\dotsc,n$,
$b\leq c$, can be written in terms of the $\frac{1}{2}n^2(n+1)-n^2
=n\binom{n}{2}$ remaining unknowns.

We can also compute the first prolongation $\mathcal{P}^1$ of \eqref{E[j,l]}.
%
Furthemore, we must impose the vanishing of the curvature, i.e.,
the equations (C) in Remark \ref {remark9}.

Therefore, by joining the systems \eqref{E[j,l]}, its first
prolongation
and the independent components of the curvature,
a system is obtained,
which has $n^2+\frac{1}{2}n^2(n+1)+\frac{1}{3}n^2(n^2-1)
=\frac{1}{6}n^2(3n+2n^2+7)$ equations in the
$\frac{1}{2}n^2(n+1)+\frac{1}{2}n^3(n+1)=\frac{1}{2}n^2(n+1)^2$
unknowns $\Gamma ^h_{ij}$,
$\frac{\partial \Gamma ^h_{ij}}{\partial x^k}$,
$1\leq i\leq j\leq n$, $h,k=1,\dotsc,n$.
We cannot write down the formulas that solve this system
because they are very involved; for example, for $n=3$,
$9$ unknown Gammas $\Gamma ^h_{ij}$, $1\leq h\leq 3$,
$1\leq i\leq j\leq 3$, can be written as functions
of the $9$ remaining Gammas, $F_j$,
and $\frac{\partial F_j}{\partial x^l}$,
$j,l=1,\dotsc,3$ in the system \eqref{E[j,l]},
whereas in the system
$\mathcal{P}^1\cup \{ R^i_{jkl}=0:1\leq k\leq j\leq n, 1\leq k<l\leq n\} $,
$36$ unknowns $\frac{\partial \Gamma ^h_{ij}}{\partial x^k}$, $h,k=1,2,3$,
$1\leq i\leq j\leq 3$, can be written as functions
of the $15$ remaining derivatives of Gamma's, $F_j$,
$\frac{\partial F_l}{\partial x^l}$,
and $\frac{\partial ^2F_l}{\partial x^k\partial x^l}$,
$j,k,l=1,2,3$, $k\leq l$, but in addition there appear $3$ constraints that
must be fulfilled identically $C_1=0$, $C_2=0$, $C_3=0$. They can be written
as functions of the coefficient
$C=-\frac{1}{2}dx^1\wedge dx^2\wedge dx^3([V_2,V_2])$, namely
\[
\begin{array}
[c]{rl}
C_1= & F_2\tfrac{\partial}{\partial x^3}\left( \tfrac{C}{F_2}\right) ,\\
C_2= & F_2\tfrac{\partial}{\partial x^2}\left( \tfrac{C}{F_2}\right) ,\\
C_3= & -\tfrac{F_1^2F_2}{F_3^2}  \left(F_1 \Gamma _{33}^1
-F_2\Gamma _{33}^2+F_3\Gamma _{33}^3\right)C
+\left(F_3\Gamma _{12}^3+F_2\Gamma _{13}^3\right) C\\
& -F_1C_2
+F_2^2\tfrac{\partial }{\partial x^1}\left(\frac{C}{F_2}\right)
-\tfrac{F_1}{F_3^2}\left( F_2\frac{\partial F_3}{\partial x^3}
-F_3\frac{\partial F_3}{\partial x^2}
-F_3\frac{\partial F_2}{\partial x^3}\right) C,
\end{array}
\]
and we can conclude that the constraints vanish
by simply applying Corollary \ref{cor}.
\end{remark}
\section{Some applications}
\begin{enumerate}
\item
Two $2$-forms $\omega _2$, $\omega _2^\prime $ on a manifold $M$ have
constant coefficients if there exist linear connections $\nabla $,
$\nabla ^\prime $ with vanishing torsion and curvature tensors such that
$\nabla \omega _2=0$, $\nabla ^\prime \omega _2^\prime =0$. Furthermore, to
say that $\omega _2$, $\omega _2^\prime $ have constant coefficients
``simultaneously'' means that $\nabla =\nabla ^\prime$.
As $\nabla $ parallelizes $\omega _2$, we have
$X(\omega _2(Y,Z))=\omega _2(\nabla _XY,Z)+\omega _2(Y,\nabla _XZ)$.
If $\omega _2$ is of maximum rank, i.e., $\operatorname*{rank}\omega _2
=2n=\dim M$, then we can define an endomorphism $J\colon TM\to TM$
by setting $\omega _2(X,J(Y))=\omega _2^\prime (X,Y)$,
$\forall X,Y\in T_xM$, and as a calculation shows, we have
$\left( \nabla _X\omega _2^\prime \right) (Y,Z)
=\omega _2\left( Y,(\nabla _XJ)Z\right) $,
$\forall X,Y,Z\in \mathfrak{X}(M)$.

Hence $\nabla $ also parallelizes $\omega _2^\prime $ if and only if it
parallelizes $J$, i.e., $\nabla J=0$. The previous result can be applied to
the study of complex valued $2$-forms $\omega _2+i\omega _2^\prime
\in \Omega ^2(M)\otimes_{\mathbb{R}}\mathbb{C}$ with constant coefficients.

\item
An almost Hermitian manifold $(M,g,J)$ is Kaehler if and only if its
fundamental $2$-form $\Phi $ is closed. The manifold $M$ is Kaehler if and only
if $d\Phi =0$. As $\Phi $ is defined by $\Phi (X,Y)=g(X,JY)$ it follows that
$\Phi $ is of maximum rank, and we can conclude from Example 1.7 above.
Moreover, if there exists an analytic coordinate system $z=(z^j)_{j=1}^n$,
$n=\dim M$, such that $\Phi =\sum_{\alpha=1}^ndx^\alpha \wedge dy^\alpha
=-\frac{1}{2i}\sum_{\alpha=1}^{n}dz^\alpha \wedge d\bar{z}^\alpha $,
with $x^\alpha =\operatorname{Re}z^\alpha $,
$y^\alpha =\operatorname{Im}z^\alpha $,
$1\leq \alpha \leq n$, then $g$ is flat, as
$\Phi =i\sum _{\alpha,\beta =1}^ng_{\alpha \bar{\beta}}dz^\alpha
\wedge d\bar{z}^\beta =$ $-\frac{1}{2i}\sum _{\alpha =1}^{n}dz^\alpha
\wedge d\bar{z}^\alpha $. Hence the
coefficients $g_{\alpha \bar{\beta}}$ are constant.

\item
Let $\omega_{p+q}$ be a form of degree $p+q$ and type $(p,q)$ on a
complex manifold $M$. If $z=(z^j)_{j=1}^n$, $n=\dim_{\mathbb{C}}M$, is an
analytic coordinate system on an open domain $U\subset M$ such that
\[
\omega_{p+q}={\textstyle\sum\limits_{i_1<_{\ldots}<i_p,j_1<_{\ldots}<j_p}}
\lambda _{i_1,_{\ldots},i_p,j_1,\dotsc,j_q}(z,\bar{z})dz^{i_1}
\wedge \ldots\wedge dz^{i_p}\wedge d\bar{z}^{j_1}\wedge \ldots \wedge
d\bar{z}^{j_q},
\]
with $\lambda _{i_1,\dotsc,i_p,j_1,\dotsc,j_q}\in \mathbb{C}$, then
the real and imaginary parts of $\omega _{p+q}$\ have constant coefficients,
as writing  $x^j=\operatorname{Re}z$,
$y^j=\operatorname{Im}z^j$, $1\leq j\leq n$, we deduce
$dz^{1}\wedge\ldots \wedge dz^{r}
=\sum _{|I|+|J|=r,I\cap J=\emptyset }\varepsilon _{IJ}(dx)^I\wedge(dy)^J$,
where $I$ and $J$ are multi-indices $I=(i_1,\dotsc,i_k)\in \mathbb{N}^k$,
$0\leq k\leq r$, $J=(j_1,\dotsc,j_l)\in \mathbb{N}^l$, $k+l=r$,
and $\varepsilon _{IJ}\in\{\pm 1,\pm i\} $. Hence
$\operatorname{Re}\omega_{p+q}$ and
$\operatorname{Im}\omega_{p+q}$ are constant coefficients when written in the
real coordinate system $(x^j,y^j)_{j=1}^n$. The converse however is not
true as shows the case of the fundamental $2$-form $\Phi $ (which is the type
$(1,1)$) of a Kaehler manifold with non-flat Riemannian metric as proved in
the previous item.

\item
If $(\omega ^1,\dotsc,\omega ^n)$ is basis of Maurer-Cartan forms on a
connected Lie group $G$, then $G$ is Abelian if and only if each of these
forms $\omega ^i$, $1\leq i\leq n$, have constant coefficients. In fact, let
$(X_1,\dotsc,X_n)$ be the dual basis to $(\omega ^1,\dotsc,\omega ^n)$
of left-invariant vector fields. If the forms $\omega ^i$, $1\leq i\leq n$,
have constant coefficients, then they are closed,
hence $0=(d\omega ^i)(X_j,X_k)=-\omega ^i([X_j,X_k])$,
$\forall i,j,k=1,\dotsc,n$. Therefore $[X_j,X_k]=0$,
thus proving that the Lie algebra of $G$ is Abelian and consequently
so is $G$, because it is connected. Conversely, if $G$ is Abelian,
then the torsion and curvature tensors of the linear
connection $\nabla ^{\mathbf{f}}$ parallelizing the linear frame
$\mathbf{f}=(X_1,\dotsc,X_n)$ vanish; hence
$\nabla ^{\mathbf{f}}=\nabla ^{\mathbf{u}}$
for some coordinate system $\mathbf{u}=(u^i)_{i=1}^n$, and we have
$\nabla ^{\mathbf{f}}\omega ^i=0$, $1\leq i\leq n$.

\item
Let $V=\mathbb{R}^n$, let
$t\in T_p^q(V)=(\otimes ^pV^\ast)\otimes (\otimes^qV)$
be a given tensor and let $G\subseteq GL(n,\mathbb{R})$
be the isotropy subgroup of $t$; namely,
$G=\{A\in GL(n,\mathbb{R}):A\cdot t=t\} $, where the dot
denotes the natural action of $GL(n,\mathbb{R})$ on
$T_p^q(V)$. If $M$ is an $n$-dimensional $C^\infty$ manifold and
$\pi \colon P\to M$ is a $G$-structure, then a tensor field $\tau $
of type $(p,q)$ on $M$ can be defined as follows: A linear frame
$u\in \pi ^{-1}(x)$ determines an isomorphism $u\colon V\to T_xM$,
which induces a $GL(n,\mathbb{R})$-equivariant isomorphism
$(u_q^p)_x\colon T_p^qV\to T_p^q(T_xM)$, and we set
$\tau (x)=(u_q^p)_x(t)$. This formula does not depend
on $u$ because of the definition of $G$. Then $P$ is
integrable if and only if each point of $M$ has a coordinate neigbourhood
$(U;\,x^1,\dotsc,x^n)$ with respect to which the components of $\tau$\ are
constant functions on $U$; e.g., see \cite{Koba}.

\item
There is a close relationship between the group of symmetries of a $p$-form
and the set of torsion-free flat linear connections that parallelize it.
For every $\omega _p\in \Omega ^p(M)$, we set
$\operatorname*{Sym}(\omega _p)=\{ \phi\in\operatorname*{Diff}M:
\phi ^\ast \omega _p=\omega _p\} $ and let $P(\omega _p)$ be
the set of coordinate systems whose associated linear
connection parallelizes $\omega _p$; namely
$P(\omega _p)=\{ \mathbf{u}=(u^1,\dotsc,u^n):
\nabla ^{\mathbf{u}}\omega _p=0\} $.

To say that $\omega _p$ has constant coefficients is equivalent to saying
that $P(\omega _p)$ is non-empty.
When making a linear change of coordinates $u^i= \alpha _j^ix^j$,
the matrix $A= (\alpha _j^i)_{i,j=1}^n$ in $GL(n,\mathbb{R})$,
a form with constant coefficients
$\omega _p$
still has constant coefficients in the new coordinate system.

Consequently, $GL(n,\mathbb{R})$ acts on the left on $P(\omega_p)$ by
composition, namely $(A,\mathbf{x})
\mapsto A\cdot\mathbf{x}=A\circ (x^1,\dotsc,x^n)$,
$\forall A\in GL(n,\mathbb{R})$,
$\forall\mathbf{x}\in P(\omega _p)$.

Once a point $x_0\in M$ and a coordinate system $\mathbf{x}=(x^1,\dotsc,x^n)$
have been fixed, every coordinate system $\mathbf{u}=(u^1,\dotsc,u^n)$
defined on an open neighbourhood $U$ of $x_{0}$ can be
written as follows: $\mathbf{u}=A\circ \mathbf{u}^\prime $, where the matrix
$A$ is as above, and $\alpha _j^i=\frac{\partial u^i}{\partial x^j}(x_0)$,
$\frac{\partial u^{\prime i}}{\partial x^j}(x_0)
=\delta _j^i, i,j=1,\dotsc,n$.
Therefore, the quotient set $P(\omega _p)/GL(n,\mathbb{R})$
can be identified to the subset $P^0(\omega _p)\subset P(\omega _p)$
of the coordinate systems whose Jacobian matrix at $x_0$
with respect to $\mathbf{x}$\ is the identity map.

In summary: Once a point $x_0\in M$ and a coordinate system $\mathbf{x}$
have been fixed, there exists a one-to-one
correspondence between the group $\operatorname*{Sym}(\omega _p)$
and the set $P^0(\omega _p)$.

The treatment for a $q$-vector field is entirely analogous.
\end{enumerate}

\noindent\textbf{Authors' addresses}

\smallskip

\noindent(J.M.M.) \textsc{Instituto de Tecnolog\'{\i}as F\'{\i}sicas y de la
Informaci\'on, CSIC, C/ Serrano 144, 28006-Madrid, Spain.}

\noindent\emph{E-mail:\/} \texttt{jaime@iec.csic.es}

\medskip

\noindent(L.M.P.C.) \textsc{Departamento de Matem\'atica Aplicada a las
T.I.C., E.T.S.I. Sistemas Inform\'aticos, UPM,  C/ Alan Turing, s/n 28031-Madrid,
Spain}

\noindent\emph{E-mail:\/} \texttt{lm.pozo@upm.es}

\medskip

\noindent(M.E.R.M.) \textsc{Departamento de Matem\'atica Aplicada, Escuela
T\'ec\-nica Superior de Arquitectura, UPM, Avda.\ Juan de Herrera 4,
28040-Madrid, Spain}

\noindent\emph{E-mail:\/} \texttt{eugenia.rosado@upm.es}
\end{document}